\input amstex.tex

\input amsppt.sty

\TagsAsMath

\magnification=1200

\hsize=5.0in\vsize=7.0in

\hoffset=0.2in\voffset=0cm

\nonstopmode

\document

\input amstex.tex
\input amsppt.sty
\TagsAsMath \NoRunningHeads \magnification=1200
\hsize=5.0in\vsize=7.0in \hoffset=0.2in\voffset=0cm \nonstopmode

\document

\topmatter
\title{$L^p$ continuity of wave operators    in $\Bbb Z$ }
\endtitle

\author
Scipio Cuccagna
\endauthor

\address
DISMI University of Modena and Reggio Emilia, via Amendola 2,
Padiglione Morselli, Reggio Emilia 42100 Italy\endaddress \email
cuccagna.scipio\@unimore.it \endemail

\abstract We recover for discrete  Schr\"odinger operators on the
lattice $\Bbb Z$, stronger analogues of the results  by Weder
\cite{W1} and by D'Ancona \& Fanelli \cite{DF} on $\Bbb R$.
\endabstract

\endtopmatter

 \head \S 1 Introduction \endhead
 We
consider  the discrete Schr\"odinger operator
$$(Hu)(n)=-( \Delta  u)(n)+q(n)u(n) \tag 1.1$$
with the   discrete Laplacian $\Delta$ in $\Bbb Z$, $ (\Delta u)(n)=
u(n+1)+u(n-1)-2u(n)  $  and a potential $q=\{ q(n), n\in \Bbb Z \}$
with $q(n)\in \Bbb R$ for all $n$. In $\ell ^2 (\Bbb Z)$   the
spectrum is $\sigma (-\Delta
 )=[0,4 ]$.
  Let for $\langle n \rangle=\sqrt{1+n^2}$
$$\aligned & \ell ^{p,\sigma}=\ell ^{p,\sigma}(\Bbb Z)=\{ u=\{ u_n\} :
\| u\|  _{\ell ^{p,\sigma}} ^p=\sum _{n\in \Bbb Z} \langle n \rangle
^{p\sigma }|u(n)|^p<\infty \} \text{ for $p\in [1,\infty )$}\\& \ell
^{\infty ,\sigma} =\ell ^{\infty ,\sigma}(\Bbb Z)=\{ u=\{ u(n)\} :
\| u\|  _{\ell ^{\infty,\sigma}}
 =\sup _{n\in \Bbb Z} \langle n \rangle ^{ \sigma }|u(n)| <\infty
\}  .\endaligned $$   We   set $\ell ^{p  } =\ell ^{p,0}.$ If $q\in
\ell ^{1,1}$ then $H$ has at most finitely many eigenvalues, see the
Appendix. The  eigenvalues are simple and are not contained  in
$[0,4]$, see for instance Lemma 5.3 \cite{CT}.  We denote by
$P_c(H)$ the orthogonal projection in $\ell ^2$ on the space
orthogonal to the space generated by the eigenvectors of $H$.
  $P_c(H)$ defines a projection in $\ell ^p$  for any
   $p\in [1, \infty ]$, see Lemma 2.6 below. We set
    $\ell ^p_c(H):=P_c(H)\ell ^p$.
By $q\in \ell ^1$ ,  $q$ is a trace class operator. Then, by
Pearson's Theorem, see Theorem XI.7\cite{RS}, the following two
limits exist in $\ell ^2$, for $w\in \ell^2_c(H )$ and $u\in \ell^2
$:

$$W u= \lim _{t\to +\infty }   e^{  it H
 }e^{  it\Delta }u \, , \quad    Z w= \lim _{t\to +\infty }  e^{ -it\Delta } e^{ -it H
 }w. \tag 1.2$$
The operators $W$ and $Z$ intertwine $-\Delta $ acting in $\ell _2$
with $H$ acting in $\ell ^2_c(H)$. Our main result  is the
following:

\proclaim{Theorem 1.1} Consider the operators $W$ initially defined
in $\ell ^2\cap \ell ^p$ and $Z$ initially defined in $\ell
^2(H)\cap \ell ^p$.

{\item {(1)}} Assume $H$ does not have resonances in 0 and 4. Then
for $q\in \ell ^{1,1}$ the operators extend into isomorphisms
$W:\ell ^p \to \ell ^p_c(H)$ and $Z:\ell ^p_c(H) \to \ell ^p   $ for
all $1<p<\infty .$ {\item {(2)}} Assume $H$ has resonances in 0
and/or  4. Then the above conclusion is true for  $q\in \ell
^{1,2}$. {\item {(3)}} Assume that $q\in \ell ^{1,2+\sigma }$ with
$\sigma
>0$. Then
$W$ and $Z$ extend into isomorphisms also for $p=1,\infty  $ exactly
when   both 0 and  4 are resonances and  the transmission
coefficient $T(\theta )$, defined for $\theta \in \Bbb T=\Bbb R/2\pi
\Bbb Z$, satisfies $T(0)=T(\pi )=1$.
\endproclaim
{\it Remark 1.}   $W$ extends into a bounded operator for
$p=1,\infty$ when the sum of the operators (3.1)--(3.4) is bounded
and this can happen only for $T(0)=T(\pi )=1$.
\medskip {\it
Remark 2.} We do not know if Claim 3 holds
  with $\sigma =0.$
\medskip
{\it Remark 3.} $\lambda =0$ or  $\lambda =4$ is a resonance exactly
if $Hu= \lambda u$ admits a nonzero solution in $\ell ^\infty $. We
say that $H$ is generic if both 0 and  4 are not resonances.

\medskip

{\it Remark 4.} Since $Z=W^*$, by duality it will be enough to
consider   $W$.
\bigskip

 Theorem 1.1   provides dispersive
estimates for solutions of the Klein Gordon equation $u_{tt}+H u
+m^2u=0$. In particular in the case of Claim 3, we obtain the
optimal $\ell^1\to \ell ^\infty $ estimate, thanks also to \cite{SK}
which deals with the $H=-\Delta$ case. The result for $T(0)=1$  by
\cite{W1} proved crucial to us
  for a nonlinear problem  in \cite{C}.
There is a close analogy between the theories in $\Bbb Z$ and in
$\Bbb R$. Claims 1 and 2  in  Theorem 1.1 are analogous to the
result in \cite{DF} for $\Bbb R$ while claim 3 is related to
analysis in \cite{W1}. Our proof mixes the approach in \cite{W1}
with estimates \cite{CT}, which in turn is inspired by
  \cite{GS,DT}. Some effort is spent proving   formulas for which
we do not know references in the discrete case. The main theme here
and in \cite{CT}, is that cases $\Bbb Z$  and $\Bbb R$ are very
similar. In particular one can see in  \cite{CT}   a theory of Jost
functions in $\Bbb Z$ very similar to the one for $\Bbb R$,
following the treatment in \cite{DT}. The present paper is inspired
by various recent  papers on dispersion theory for the group
$e^{itH} $, see \cite{SK,KKK,PS,CT}. In particular the
 bound $| e^{it\Delta }(n,m)|\le C \langle t \rangle ^{-1/3}$ was
proved   in \cite{SK}.   The bound $|P_c(H)e^{itH}(n,m)|\le C
\langle t \rangle ^{-1/3}$ was proved in \cite{PS} for $q\in \ell
^{1,\sigma }(\Bbb Z)$ with $\sigma
>4$ and for $H$ without  resonances. This result
was extended by \cite{CT} to $q\in \ell ^{1,1}$ for $H$ without
resonances and to $q\in \ell ^{1,2}$  if 0 or 4 is a resonance.
\cite{CT} is able  produce for $\Bbb Z $ essentially the same
argument introduced in \cite{GS} for $\Bbb R $, thanks to a a theory
of Jost functions in $\Bbb Z $  which is basically the same of that
for $\Bbb R $. Here we recall that \cite{GS} for Schr\"odinger
operators on $\Bbb R$ improves  an earlier result in \cite{W2}.
Theorem 1.1 is the natural transposition to  $\Bbb Z $, with some
improvements, of the theory of wave operators for $\Bbb R$ in
\cite{W1,GY,DF}.  We simplify the argument in \cite{DF} for claims
(1) and (2) of Theorem 1.1 and, for claim (3), we use weaker decay
hypotheses on the potential than \cite{W1}.

\bigskip
We end with some notation. Given an operator $A$ we set
$R_A(z)=(A-z)^{-1}$.     $\Cal S(\Bbb Z)$  is the set of functions
$f:\Bbb Z \to \Bbb R$ with $f(n)$
 rapidly decreasing as $|n|\nearrow \infty$. For $u\in \ell ^2$ we
 set
$F_0[u](\theta ):=\frac{1}{\sqrt{2\pi }}\sum _{n\in \Bbb Z}
e^{-in\theta} u(n).$ We set $\Bbb T=\Bbb R/2\pi \Bbb Z$. $2\Bbb Z$
is the set  of even integers; $2\Bbb Z+1$ is the set  of odd
integers.  We set
$$\eta (\mu )=\sum _{\nu =\mu }^{\infty} |q(\nu )|  \text{ and }
\gamma (\mu )=  \sum _{\nu =\mu }^\infty (\nu -\mu ) |q(\nu )| .$$
Given $f\in L^1(\Bbb T)$ we set $\widehat{f}(\nu )=\int _{-\pi}^\pi
e^{-i\nu \theta}f(\theta ) d\sigma $, with $d\sigma =d\theta
/\sqrt{2\pi}.$

\bigskip
 \head \S 2 Fourier transform associated to $H$ \endhead

We recall that the resolvent $R_{-\Delta}(z)$  for $ z\in \Bbb C
\backslash [0,4]$ has kernel
$$
R_{-\Delta}(m, n, z)=\frac{-i}{2\sin \theta }e^{-i \theta   |n-m|},
\ \ m,n \in \Bbb Z,
$$
with $\theta$ a solution     to $    2(1-\cos \theta ) =z
 $ in
$   D=\{ \theta : -\pi \le \Re \theta \le \pi , \, \Im \theta <0 \}.
$  In \cite{CT}  it is detailed the existence of functions ${f}_{
\pm } (n, \theta )$ with
$$H{f}_{  \pm } (\mu , \theta )=z{f}_{  \pm } (\mu, \theta ) \text{ with }
\lim _{\mu \to \pm \infty}\left [ {f}_{  \pm } (\mu , \theta ) -
e^{\mp i \mu \theta }\right  ]=0.\tag 2.1$$ We have
$$\aligned   {f}_{  \pm } ( \mu , \theta )=&
 e^{\mp in\theta }  - \sum _{\nu =\mu }^{\pm \infty
} \frac{\sin  (\theta (\mu -\nu ))}{\sin \theta } q(\nu )  {f}_{ \pm
} ( \nu , \theta )  .\endaligned \tag 2.2
$$
Define  $
 {m}_{  \pm} $    by  $ {f}_{  \pm } (n , \theta )=
 e^{\mp i n\theta
  } {m}_{  \pm } (n,\theta ).$
Lemma  5.1  \cite{CT}   implies that for   fixed $n$
$$m_{\pm }(n,\theta )=1+\sum _{\nu =1}^{\infty} B_{\pm }(n ,\nu )
e^{-i\nu \theta }.\tag 2.3$$ In Lemma 5.2  \cite{CT} it is proved:

\proclaim{Lemma 2.1} For $q\in \ell ^{1,1}$ and setting
 $ B _{+}(n ,0 )=0$ for all $n$,
we have
$$\aligned & B _{+}(n  ,2\nu  ) =\sum _{l=0}^{\nu -1 }
\sum _{j=n+\nu -l}^{\infty} q(j) B _{+}(j,2l+1  )\\& B _{+}(n  ,2\nu
-1 ) =\sum _{l=n+\nu}^{ \infty}q(l)  +\sum _{l=0}^{\nu -1 }\sum
_{j=n+\nu -l}^{\infty} q(j) B _{+}(j,2l   ).
\endaligned
$$ We have for $n\ge 0$  the estimate
  $|  B_+(n, \nu)| \leq \chi _{[1,\infty )}( \nu)
  e^{\gamma(0)}\eta( \nu) $.
Similarly for $n\le 0$  we have
  $|  B_-(n, \nu)| \leq \chi _{[1,\infty )}( \nu)
  e^{\widetilde{\gamma}(0)}\widetilde{\eta}( \nu) $
with $\widetilde{\gamma}(\mu)$ and $\widetilde{\eta}( \mu) $ defined
like $ {\gamma}(\mu)$ and $ {\eta}( \mu) $ but with $q(\nu )$
replaced by $q(-\nu )$.
\endproclaim

Lemma 2.1 implies what follows, see the proof of Lemma 5.10
\cite{CT}: \proclaim{Lemma 2.2}  If $q\in \ell ^{1,1+\sigma }$ for
$\sigma \ge 0$, then
  $\|  B_\pm (n, \cdot )\| _{\ell ^{1,\sigma }}
  \le C_\sigma  \|  q\| _{\ell ^{1,1+\sigma }}  $ for $\pm n \ge 0.$
\endproclaim
We recall that for   two given functions $u(n)$ and $v(n)$ their
Wronskian is $[u,v](n)= u(n+1)v(n)-u(n) v(n+1)$. If $u$ and $v$ are
  solutions of $Hw=zw$    then $[u,v]$ is constant.
In particular we set $W(\theta ):=[
 f_+(\theta ), f_-(\theta )]$ and $W_1(\theta ):=[
 f_+(\theta ), \overline{f}_-(\theta )]$. By an argument in Lemma
 5.10 \cite{CT} we have:
\proclaim{Lemma 2.3} If for $\sigma \ge 0$ we have $q\in \ell
^{1,1+\sigma }$,  then
  $W(\theta ), W_1(\theta )\in \ell ^{1,\sigma }$.
\endproclaim
 Lemma 5.4
  \cite{CT} states:

\proclaim{Lemma 2.4} Let $q\in \ell ^{1,1}$. For $\theta \in [-\pi ,
\pi ]$ we have $\overline{f_{\pm}(n,\theta )} =f_{\pm}(n,-\theta )$
and for $\theta \neq 0,\pm \pi $ we have
$$ {f}_{  \mp } (n,\theta )=
\frac{1}{ T  (\theta )}\overline{ {f}_{  \pm } (n,\theta )}
+\frac{R_{ \pm} (\theta)}{T  (\theta  )}  {f}_{  \pm } (n,\theta )
\tag 1 $$ where $T  (\theta )$ and $R_\pm (\theta  )$ are defined by
(1) and satisfy: $$  \align   & [\overline{ {f}_{ \pm } (\theta )},
{f}_{ \pm } (\theta )] = \pm 2i\sin \theta , \tag 2\\& T  (\theta
)= \frac {- 2 i\sin \theta } {W(\theta )  } \, , \quad R_{  + }
(\theta ) = -\frac {\overline{W}_1(\theta )  } {W(\theta ) }\, ,
\quad R_{ + } (\theta ) = -\frac { {W}_1(\theta )  } {W(\theta )
}\tag 3
\\&
  \overline{T  (\theta  )}=T  (-\theta  ) \, , \,
\overline{R_{  \pm } (\theta  )}=R_{  \pm } (-\theta  ) , \tag 4
\\&|T  (\theta  )|^2+ |R_{ \pm } (\theta ) |^2=1\, ,
\quad T (\theta )  \overline{R_{ \pm } (\theta  )}+R_{  \mp} (\theta
) \overline{ T  (\theta  )}=0.\tag 5\endalign
$$
\endproclaim
Lemma 5.5 \cite{CT} states:

 \proclaim{Lemma 2.5}

 {\item {(1)}}   For $\theta \in [-\pi ,\pi ]\backslash
 \{ 0, \pm \pi \}$ we have $W(\theta )\neq 0$.
 We have $|W(\theta )|\ge 2|\sin \theta |$ for all
 $\theta \in [-\pi , \pi ]$
 and  in the
generic case $|W(\theta )|>0$.

{\item {(2)}} For  $j=0,1$ and  $q\in \ell ^{1,1 +j}$ then $W(\theta
)$ and $W_1(\theta )$  are in  $  C^j[-\pi , \pi ]$.

{\item {(3)}}  If $q\in \ell ^{1,2}$ and $ W(\theta _{0})=0$ for a $
\theta _{0}\in \{ 0, \pm \pi \}$, then $\dot W(\theta _{0})\neq 0$.
In particular if $q\in \ell ^{1,2}$, then $T  (\theta  )=  - { 2
i\sin \theta }/ {W(\theta )  }$  can be extended continuously in
$\Bbb T$.

\endproclaim
We have the following result: \proclaim{Lemma 2.6} Assume that $q\in
\ell ^{1,1 }$ if $H$ is generic and $q\in \ell ^{1,2  }$ if $H$ has
a resonance at 0 or at 4. Then the following statements hold.

{\item {(1)}}
 $H$ has finitely many eigenvalues.
{\item {(2)}} If $\lambda$ is an eigenvalue, then $\text{dim} \ker
(H-\lambda )=1$. {\item {(3)}}
  If there are eigenvalues
they are in $\Bbb R \backslash [0,4]$. {\item {(4)}} Let $\lambda
_1$,...,$\lambda _n$ be the eigenvalues and $\varphi
_1$,...,$\varphi _n$ corresponding eigenvectors with $\| \varphi
_j\| _{\ell ^2}=1$. Then for fixed $C>0$ and $a>0$ we have $|
\varphi _j (\nu )|\le C e^{-a|\nu|}$ for all $j=1,...,n$ and for all
$\nu \in \Bbb Z$. {\item {(5)}} Let $P_d(H):=\sum _j\varphi _j
\langle \quad , \varphi _j\rangle $. Then $P_d(H)$ and $
P_c(H):=1-P_d(H)$ are bounded operators in $\ell ^p$ for all $p \in
[1,\infty ].$
\endproclaim
{\it Proof.} (1) is proved in the Appendix. (2) and (3) are in Lemma
5.3 \cite{CT}. (5) follows from (4). (4) follows from the fact that
by the proof in Lemma 5.3 \cite{CT} there are constants $A(\pm , j)$
such that $\varphi _j (\nu )=A(\pm , j) f_{\pm }(\nu , \theta _j)$,
with $\theta _j\in D$ such that $\lambda _j=2(1-\cos (\theta _j))$.
The fact that $\lambda _j\not \in [0,4]$ implies $\Im (\theta _j)<0$
for all $j$.

\bigskip
By Lemmas 5.6-9 \cite{CT}  we have
$$\aligned & P_c( H )u=
 \frac 1{2\pi i}  \int _{ 0}^{4}
 \left [ R_{H  }^+(\lambda
)- R_{H  }^-(\lambda   )
 \right ] ud\lambda =\\& =\frac{1}{2\pi i}
\sum _{\nu \in \Bbb Z}\int _{ -\pi  }^\pi
    K  (n, \nu ,\theta ) d \theta u(\nu ) \text{ with }
\endaligned \tag 2.4
$$

$$ \aligned & K  (n, \nu ,\theta )=
 {f_-}  (n, \theta )   f_+  (\nu , \theta )  \frac {\sin ( \theta
) } {W(\theta )} \text{ for $\nu >n$}\\& K  (n, \nu ,\theta )=
 {f_+}  (n, \theta )   f_-  (\nu , \theta )  \frac {\sin ( \theta
) } {W(\theta )} \text{ for $\nu \le n$.}
\endaligned \tag 2.5$$
Consider now plane waves defined as follows:

 \proclaim{Definition
2.7 } We consider the following functions: $$\aligned & \psi (\nu
,\theta )=\frac{1}{\sqrt{2\pi }}T(  \theta ) e^{-i\nu \theta } m_+
(\nu, \theta )\text{ for $  \theta \ge 0$ } \\& \psi (\nu ,\theta
)=\frac{1}{\sqrt{2\pi }}T(- \theta ) e^{-i\nu \theta } m_- (\nu , -
\theta )\text{ for $  \theta < 0$ } .\endaligned $$
\endproclaim

 \proclaim{Lemma
2.8} The kernel  $P_c(H)(\mu ,\nu)$  of $P_c(H)$ can be expressed as
$$P_c(H)(\mu ,\nu)=\int _{-\pi}^{\pi} \overline{\psi (\mu
,\theta )}  \psi (\nu ,\theta ) d\theta .\tag 1$$
\endproclaim
{\it Proof.} We assume $\mu \ge \nu$. By (2.4-5)
$$P_c(H)(\mu ,\nu)=\frac{1}{2\pi i}\int _{0}^{\pi}\left [
   \frac {{f_-}  (\nu , \theta )
  f_+  (\mu , \theta )} {W(\theta )}-  \frac
{{f_-}
   (\nu ,- \theta )
    f_+  (\mu ,- \theta ) } {W(-\theta )} \right ]
     \sin ( \theta )  \, d\theta .$$
We have by Lemma 2.4

$$\aligned &  \overline{f_{
\pm } (n,\theta )}= f_{ \pm } ( n,-\theta)\, , \, \overline{T(\theta
)}=T(-\theta ) \, , \, \overline{R_{\pm } (\theta )}=R_{\pm }
(-\theta),
\\&    f_{ - } (\nu , -\theta )  =   {T  (\theta)
 f_{ + } (\nu , \theta )
-  R_{- }  (\theta ) f_{ -} (\nu ,\theta )}  ,\\& f_{ + } (\mu ,
  \theta)=   \overline {T (\theta )  {f_{ - } (\mu ,\theta )} -R_{+ } ( \theta
) {f_{ +} (\mu ,\theta )}} .\endaligned
$$
Substituting the last two lines in the square bracket  in the
integral,

$$\align   &  [\cdots ]=    \frac{  \overline {T(\theta )
   f_{ - }  (\mu ,\theta )}     f_{  -}(\nu , \theta ) }
  {W( \theta )} - \frac{   T(\theta )
     f_{ +}  (\nu ,  \theta)   f_{  +} (\mu , -\theta ) }
   {W( -\theta )}\tag 2 \\&
- \overline{f_{ +} (\mu ,  \theta )}   f_{  -}(\nu ,   \theta )
\left [ \frac{ \overline{R_{+}( \theta)}  } {W(  \theta ) } - \frac{
R_{-}(\theta ) } {W(- \theta )} \right ].
\endalign  $$
The last line is  zero by  (5) Lemma 2.4 and by
$$ -i\sin (\theta )\left [ \frac{ \overline{R_{+}( \theta)}  } {W(  \theta ) } - \frac{
R_{-}(\theta ) } {W(- \theta )} \right ]= (T \overline {R_+}
+\overline{T} R_-  )(\theta) =0.$$ We have by $T(\theta )=-i\sin
(\theta )/W(\theta )$
$$\text{rhs}(2)= \frac{1}{2\pi } |T(\theta )|^2
\overline{ f_{ + }  (\mu ,\theta )}     f_{  +}(\nu , \theta )
+\frac{1}{2\pi } |T(\theta )|^2    \overline{f_{ - }  (\mu , \theta
)}
 f_{ -}(\nu ,  \theta ) .$$ This yields formula (1) for
$\mu \ge \nu$. For $\mu < \nu$ the argument is similar.

\proclaim{Lemma 2.9}  Let $F[u](\theta ):=\sum _n    \psi
  (n,\theta ) u(n).$ Then:  {\item {(1)}} $F:\ell ^2_c(H)\to
L^2(\Bbb T )$ is an isometric isomorphism. {\item {(2)}}
$F^{\ast}[f](n):=
 \int _{-\pi}^{\pi }  \overline{\psi
(n,\theta )} f(\theta ) d\theta  $ is the inverse of $F$. {\item
{(3)}} $F[Hu](\theta )= 2(1-\cos \theta ) F[u](\theta ).$
\endproclaim
  $F[u](\theta )$ is a generalization of Fourier series
expansions   $F[u_0](\theta )$. Lemma 2.9 is a consequence of Lemma
2.8 except for the fact  that we could have $F(\ell
^2_c(H))\subsetneqq L^2(\Bbb T  )$. The fact $F(\ell ^2_c(H))=
L^2(\Bbb T  )$ follows from $F_0(\ell ^2 )= L^2(\Bbb T  )$, from the
fact that  $W$ and $Z$ in (1.2) are isomorphisms between $\ell ^2$
and $\ell ^2_c(H)$ and from Lemma 2.10 below. In the next section
the following formula will be important:

\proclaim{Lemma 2.10} For the operator in (1.2) we have $W=F^\ast
F_0$.
\endproclaim

  We have, for $u,v\in \Cal S(\Bbb
Z)$ and $v\in L^2_c(H)$
$$\langle W u,v\rangle _{\ell ^2}- \langle  u,v\rangle _{\ell ^2}=
i\lim _{\epsilon \searrow 0} \int _{0}^{\infty }  \langle e^{itH}q
e^{it\Delta } u,v\rangle _{\ell ^2}e^{-\epsilon t} dt.$$ We have for
$L ^2=L ^2(\Bbb T)$
$$\aligned &\langle e^{itH}q
e^{it\Delta } u,v\rangle _{\ell ^2} =\langle e^{i2t(1-\cos \theta
)}F[q e^{it\Delta }  u],F[v]\rangle _{L ^2 }=\langle  F[q e^{it(
\Delta +2(1-\cos \theta )}  u],F[v]\rangle _{L ^2 }  .
\endaligned $$ Then
$$ i\int _{0}^{\infty }  \langle e^{itH}q
e^{it\Delta } u,v\rangle _{\ell ^2}e^{-\epsilon t} dt= \langle  F[q
R_{-\Delta } (2 -2\cos \theta +i\epsilon )  u],F[v]\rangle _{L ^2 }
$$ and

$$\align  & \langle W u,v\rangle _{\ell ^2}- \langle
 u,v\rangle _{\ell ^2}= \\&  =  \int _{-\pi  }^{\pi}d\theta \,
\overline{F[v]}(\theta )
 \sum _{\nu \in \Bbb Z}    \psi
(\nu,\theta )  q(\nu ) (R _{-\Delta }^+(   2  -2\cos \theta    )
u)(\nu )=\\&
 \int _{-\pi  }^{\pi} {d\theta} \,
\overline{F[v]}(\theta )\sum _{  \nu ' \in \Bbb Z}
 u(\nu ' ) \frac {-i}{2 \sin |\theta |}\sum _{\nu     \in \Bbb Z}
   e^{-i|\theta |\, \, |\nu  -\nu '  |}q(\nu ) \psi
(\nu,\theta )  .\tag 1
\endalign  $$
We claim we have
$$\psi (\mu , \theta )= e^{-i\mu \theta }/\sqrt{2\pi}+\frac{i}{2\sin \theta}
\sum _{\nu     \in \Bbb Z}
   e^{-i \theta  \, \, |\nu  -\mu  |}q(\nu )\psi
(\nu,\theta ) \text{ for $\theta >0$}\tag 2$$
$$\psi (\mu , \theta )= e^{  -i\mu  \theta }/\sqrt{2\pi}-\frac{i}{2\sin \theta}
\sum _{\nu     \in \Bbb Z}
   e^{  i \theta  \, \, |\nu  -\mu  |}q(\nu )\psi
(\nu,\theta ) \text{  for $\theta <0$}.\tag 3$$ Assuming (2)--(3)

$$\aligned & \langle W u,v\rangle _{\ell ^2}- \langle
 u,v\rangle _{\ell ^2}=
 \int _{-\pi  }^{\pi}\sum _{  \nu ' \in \Bbb Z} {d\theta} \,
\overline{F[v]}(\theta )
 u(\nu ' )\left [e^{-i\nu '\theta }/\sqrt{2\pi}-
   \psi (\nu ', \theta )  \right ] \\& =  \int _{-\pi  }^{\pi}
   {d\theta} \, \overline{F[v]}(\theta ) \left [
   {F_0[u]}(\theta )-{F[u]}(\theta ) \right ]
    = \langle F^*F_0u, v\rangle _{\ell ^2}
    -\langle  u, v\rangle _{\ell ^2}.
\endaligned $$
This yields $W=F^*F_0.$ Now we focus on (2) and (3). For $\theta >0$
it is possible to rewrite (2.2) as follows, for some constant
$A(\theta )$,
$$f_+(\mu , \theta )= e^{-i\mu \theta } A(\theta )-R_{-\Delta }^+(2-2\cos
\theta )qf_+(\cdot , \theta )(\mu ).\tag 4$$
   Using (2.2) for $f_-$ we obtain
 $-2i\sin (\theta ) A(\theta ) =[f_+(  \theta ), f_-(\mu , \theta
 )]$.
Hence  $A(\theta )=1/T(\theta )$. So multiplying  (4)  by $T(\theta
)/\sqrt{2\pi}$ we obtain  (2).
 We have for
$\theta <0$
$$f_-(\mu , \theta )= e^{ i\mu \theta } B(\theta )-R_{-\Delta }^-(2-2\cos
\theta )qf_-(\cdot , \theta )(\mu )\tag 5$$
 for some constant $B(\theta )$. One checks that
 $-2i\sin (\theta ) B(\theta ) =[f_+(  \theta ), f_-(\mu , \theta
 )]$.
Hence $B(\theta )=1/T(\theta )$. So multiplying  (5)  by $T(\theta
)/\sqrt{2\pi}$ we obtain
$$\frac{T(\theta )}{\sqrt{2\pi}}f_-(\mu , \theta )=
\frac{e^{ i\mu \theta }}{\sqrt{2\pi}}  -R_{-\Delta }^-(2-2\cos
\theta )q\frac{T(\theta )}{\sqrt{2\pi}}f_-(\cdot , \theta )(\mu ).$$
Taking complex conjugate we obtain (3).

\bigskip
\head \S 3 Bounds on $W$ \endhead

 It is not restrictive to consider $\chi
_{[0,\infty ]} (n) Wu(n)$ instead of $Wu(n)$. Indeed the proof for
$\chi _{(-\infty ,0)} (n) Wu(n)$ is similar.     Claims 1 and 2 in
Theorem 1.1 are a consequences of Lemma 3.1 below. We follow
\cite{W1}, exploiting at some crucial points  results proved in
\cite{CT} and inspired by \cite{GS}.  We set  $n_\pm (\mu , \theta
):= m_\pm (\mu , \theta )-1$.

\proclaim{Lemma 3.1}  Let $q\in \ell ^{1,1}$ in the generic case and
$q\in \ell ^{1,2}$ in the non  generic case. Then $  \| \chi
_{[0,\infty ]}  Wu \| _{\ell ^p}\le C_p \|  u \| _{\ell ^p} \quad
\forall \, p\in (1,\infty ).$
\endproclaim
{\it Proof.}
 Recall $F_0^\ast [n_\pm (\mu , \cdot )] (\nu )=B_\pm
(\mu , \nu ).$ Furthermore in Lemma 5.10 \cite{CT} it is proved that
$F_0^\ast [T]\in \ell ^1$. One can prove similarly  that also
$F_0^\ast [R_\pm ]\in \ell ^1$. For $d\sigma = d\theta /\sqrt{2\pi}$
and by $\overline{m_\pm} (\mu , \theta )={m_\pm} (\mu , -\theta )$,
$\overline{T}(\theta ) = {T}(-\theta )$, we consider
$$\aligned & Wf(\mu ) = \int _{-\pi }^\pi
\overline{\psi  (\mu , \theta )} F_0[f] ( \theta ) d\theta  = \int
_{0}^\pi T(-\theta ) e^{ i\mu \theta } m_+(\mu , -\theta ) F_0[f]
(\theta ) d\sigma
\\&+ \int _{-\pi }^0 T(  \theta ) e^{  i\mu \theta } m_-(\mu ,
\theta ) F_0[f] ( \theta ) d\sigma .\endaligned
$$
We consider only $\mu \ge 0$.
 We substitute $n_\pm (\mu , \theta ):= m_\pm (\mu , \theta )-1$
 and
$  T( \theta ) {m}_{  - } (\mu , \theta ) =
  m_{  + } (\mu , -\theta )
 +e^{ -2i\mu \theta } R_{ +} ( \theta)  {m}_{  + } (\mu , \theta )
$ obtaining

$$\aligned &\chi
_{[0,\infty ]} (\mu )Wf(\mu ) = \int _{-\pi }^\pi e^{ i\mu \theta
}T(-\theta ) \frac{1+\text{sign} (\theta )}{2}  F_0[f] (\theta )
d\sigma \\&+ \int _{-\pi }^\pi
  e^{  i\mu \theta }  \frac{1-\text{sign}
   (\theta )}{2}  F_0[f] ( \theta )
d\sigma   + \int _{-\pi }^\pi  e^{ - i\mu \theta }R_+(  \theta )
 \frac{1-\text{sign}
   (\theta )}{2} F_0[f] ( \theta ) d\sigma
 \\& + \int
_{-\pi }^\pi e^{ i\mu \theta }T(-\theta )  n_+(\mu , -\theta )
\frac{1+\text{sign}(\theta )}{2}F_0[f] (\theta ) d\sigma \\& + \int
_{-\pi }^\pi
  e^{  i\mu \theta } n_+(\mu , -  \theta )
  \frac{1-\text{sign}(\theta )}{2} F_0[f] ( \theta )
d\sigma\\& +  \int _{-\pi }^\pi  e^{ - i\mu \theta } R_+(  \theta )
n_+(\mu , \theta ) \frac{1-\text{sign}(\theta )}{2} F_0[f] ( \theta
) d\sigma .\endaligned $$ We have $\chi _{[0,\infty ]} (\mu )Wf(\mu
)= \widetilde{W}_1f(\mu )+ \widetilde{W}_2f(\mu )$ where, for
$W_j=2\sqrt{2\pi }\widetilde{W}_j$ for $j=1,2$:

$$\aligned & W_1f(\mu ) = \int _{-\pi }^\pi
e^{ i\mu \theta }    T(-\theta )   F_0[f] (\theta ) d\theta
 + \sqrt{2\pi} f
  +  \int _{-\pi }^\pi  e^{-  i\mu \theta }R_+(  \theta )
    F_0[f] ( \theta ) d\theta
  \\& + \int
_{-\pi }^\pi e^{ i\mu \theta }\left ( T(-\theta ) +1\right ) n_+(\mu
, -\theta )
   F_0[f] (\theta ) d\theta     +
     \int _{-\pi }^\pi  e^{ - i\mu \theta } R_+(  \theta )
n_+(\mu , \theta )   F_0[f] ( \theta ) d\theta ;\endaligned $$

$$\aligned & W_2f(\mu ) = \int
_{-\pi }^\pi e^{ i\mu \theta }\left ( T(-\theta )-1\right )  m_+(\mu
, -\theta )
 \text{sign}(\theta ) F_0[f] (\theta ) d\theta  -\\&
  - \int _{-\pi }^\pi  e^{ - i\mu \theta }R_+(  \theta
)m_+(\mu , \theta )  \text{sign}(\theta )  F_0[f] ( \theta ) d\theta
  .\endaligned $$ $W_1$ is bounded  for $p\in [1,\infty ]  $.
Indeed for example,
$$ \aligned & \left \| \chi _{[0,\infty )} (\cdot  )
F_0^\ast \left[ R_+(  \theta ) n_+(\mu , \theta )   F_0[f] ( \theta
)\right ] (-\, \cdot ) \right \| _{\ell ^p  }\le \\& \left \| \chi
_{[0,\infty )} (\cdot ) \left ( \left | F_0^\ast \left[ R_+  \right
] \right |
*  \chi _{[1,\infty )}e^{\gamma (0)}\eta  *  |f| \right ) (-\,\cdot
) \right \| _{\ell ^p }\\& \le  e^{\gamma (0)}\gamma (0)\| F_0^\ast
\left[ R_+ \right ]\| _{\ell ^1 }   \, \| f \| _{\ell ^p },
\endaligned $$
where we have used $|B_+(\mu , \nu )|\le \chi _{[1,\infty )}(\nu )
e^{\gamma (0)}\eta (\nu )$ for $\mu \ge 0$. Other terms of $W_1$ can
be treated similarly. By the same argument  $W_2$  is bounded  for
$p\in (1,\infty) $. For $W_2$ we cannot include $p=1,\infty$ because
$\text{sign}(\theta )$ is the symbol of the Calderon-Zygmund
operator
$$\Cal H v (\nu )=\int _{-\pi}^\pi  e^{ i\nu \theta } F_0[v](\theta
) \, d\sigma =  \frac{2i}{\pi }\sum _{\nu ' \in \nu +2\Bbb
Z+1}\frac{v(\nu ')}{\nu -\nu '}$$ which is unbounded in $\ell ^1$
and in $\ell ^\infty$.   So the proof of Lemma 3.1 is completed.
\bigskip Consider now $W_2f(\mu )= \chi _{[0,\infty
]} (\mu )W_2f(\mu )$

\proclaim{Lemma  3.2} Let $q\in \ell ^{1,2+\sigma }$ with $\sigma
>0$. Then $W_2$ extends  into a bounded operator also for
$p=1,\infty $ exactly when  both  0 and  4 are resonances  and the
transmission coefficient $T(\theta )$ defined in $\Bbb T$ satisfies
$T(0)=T(\pi )=1$.
\endproclaim
{\it Proof.} We consider a partition of unity $1=\chi +(1-\chi )$ on
$\Bbb T $ with $\chi $ even,  $\chi
 =1$ near $0$ and $\chi
 =0$ near $\pi $. Correspondingly we have $W_2=U_1+U_2$ with $U_1$
 written below and $U_2$ given by the same formula with $\chi$
 replaced by $1-\chi $.
We focus on  $U_1 $. We have $U_1=U_{11}+U_{12} $ with for $\mu \ge
0$

$$\aligned & U_{11}f(\mu ) = U_{111}f(\mu )+ U_{112}f(\mu )\\&
U_{111}f(\mu )= m_+(\mu , 0 )\int _{-\pi }^\pi  e^{ i\mu \theta }
 \left ( T(-\theta ) -T(0)\right )   \text{sign}(\theta ) \chi (\theta ) F_0[f] (\theta )
  d\theta  \\& - m_+(\mu , 0 )\int _{-\pi }^\pi
  e^{ - i\mu \theta }\left ( R_+(  \theta
)-R_+(0) \right )  \text{sign}(\theta )  F_0[f] ( \theta ) d\theta
\\& U_{112}f(\mu )= \int _{-\pi }^\pi e^{ i\mu \theta } \left
(T(-\theta ) -1\right ) \left ( n_+(\mu , -\theta ) -n_+(\mu , 0
)\right ) \text{sign}(\theta ) \chi (\theta )F_0[f] (\theta )
d\theta   \\& - \int _{-\pi }^\pi e^{ - i\mu \theta }  R_+( \theta )
 \left ( n_+(\mu , \theta )- n_+(\mu , 0 )\right )  \text{sign}(\theta ) \chi (\theta
)F_0[f] ( \theta ) d\theta
\endaligned $$
and $$  \aligned & U_{12}f(\mu ) = \chi _{[0,\infty )}(\mu )\left (
T(0)-1\right ) m_+(\mu , 0 ) \int _{-\pi }^\pi
  e^{  i\mu \theta } \text{sign} (\theta )  \chi (\theta )
    F_0[f] ( \theta )
     \\&  - \chi _{[0,\infty )}(\mu ) R_+(  0 )m_+(\mu ,
0 ) \int _{-\pi }^\pi
  e^{  -i\mu \theta } \text{sign} (\theta ) \chi (\theta ) F_0[f] (  \theta )
d\theta  \\& = \chi _{[0,\infty )}(\mu )\left ( T(0)-1\right )
m_+(\mu , 0 ) (\Cal H f)(-\mu )
       - \chi _{[0,\infty )}(\mu ) R_+(  0 )m_+(\mu ,
0 ) (\Cal H f)( \mu )  .\endaligned \tag 3.1$$ We have:
\proclaim{Lemma 3.3}   $U_{12}\in B (L^p,L^p)$  for all $p\in
[1,\infty ]$ if and only if
$$T(0)-1+R_+(  0 )=0.\tag 1$$
\endproclaim
{\it Proof.}   We have $m_+(\mu , 0 )\to 1$ for $\mu \nearrow \infty
$ if $q\in \ell ^{1,1}$.  We have $(\Cal H f)(-\mu )= (\Cal H f (-\,
\cdot ))( \mu ).$ Set $\widehat{\chi}   =F_0^\ast (\chi )$. Then
$U_{12}\in B (L^p,L^p)$   for $p=1,\infty$ exactly if
$$\align  & \chi _{\Bbb N}(\mu )\left ( T(0)-1+R_+(  0 )\right )
  \Cal H (\widehat{\chi} * f)( \mu )\in \ell ^p \text{ for all
  $f$ even in $\ell ^p$}\tag 2\\&
  \chi _{\Bbb N}(\mu )\left ( T(0)-1-R_+(  0 )\right )
 \Cal H (\widehat{\chi} *f)( \mu )\in \ell ^p \text{ for all
  $f$ odd in $\ell ^p$.}\tag 3
\endalign  $$
We show that (2) requires (1). We have $\widehat{\chi} *\chi _{\{ 0
\}} =\widehat{\chi}$ and

  $$\aligned & (\Cal H \widehat{\chi})( \mu )=  \frac{2i}{\pi \mu}
   \sum _{   \nu \in \mu
+2\Bbb Z+1}
  \widehat{\chi}(\nu ) -\frac{2i}{\pi }\sum _{   \nu \in \mu
+2\Bbb Z+1} \left [ \frac{1}{\mu}-
  \frac{1}{\mu -\nu }\right ]\widehat{\chi}(\nu ).\endaligned $$
The second term on the right is in $\ell ^1([1,\infty )$ but the
first is $ i\frac{ \sqrt{2}}{\sqrt{\pi } \mu } $, which is not in
$\ell ^1([1,\infty )$. Hence we need equality  (1). So (2) requires
(1). We now show that (3) occurs always.  It is enough to
 prove $\Cal Hf\in \ell ^p$ for all $f$ odd.
We have
$$\aligned & \sum _{   \nu \in \mu +2\Bbb Z+1}
  \frac{1}{\mu -\nu }f(\nu )=2 \sum _{   \nu \in \mu +2\Bbb Z+1}
  ^{\nu >0}
  \frac{\nu }{\mu  ^2-\nu ^2 }  f(\nu )
.\endaligned $$ So
$$\aligned &\| \Cal Hf \| _{\ell ^1} \lesssim  \sum _{\nu >0}
|f(\nu )| \sum _{   \mu \in \nu +2\Bbb Z+1}
  \frac{\nu }{|\mu  ^2-\nu ^2| }\le C \|  f \| _{\ell ^1} \endaligned $$
for a fixed $C<\infty $.
\bigskip
Our next step is to show in Lemma 3.4 that $U_{111}\in B (L^p,L^p)$
for all $p\in [1,\infty ]$. In Lemma 3.5 that $U_{112}\in B
(L^p,L^p)$  for all $p\in [1,\infty ]$. Hence $U_1\in B (L^p,L^p)$
for all $p\in [1,\infty ]$ exactly if $U_{12}\in B (L^p,L^p)$  for
all $p\in [1,\infty ]$.

  \proclaim{Lemma 3.4}  Let $q\in \ell ^{1,2+\sigma }$ with
$\sigma
>0$. Then $U_{111}\in B (L^p,L^p)$  for all $p\in
[1,\infty ]$.
\endproclaim
{\it Proof.}  If    for   $g=\left ( R_+( \theta ) -R_+(0) \right
)\text{sign}(\theta ) \chi (\theta )$ and $f=\left ( T( \theta )
-T(0) \right )\text{sign}(\theta ) \chi (\theta )$ we have $ F_0^*f
$ and $  F_0^*g\in \ell ^1$, then by $|m_+(\mu , 0)|\le C$ for all
$\mu \ge 0$, we get Lemma 3.3. Here consider only $F_0^*f $ only,
since the proof for $ F_0^*g$ is similar. We have for
$\widetilde{\chi}(\theta )$ another even smooth cutoff function in
$\Bbb T$ with $\widetilde{\chi}=1$ on the support of $\chi$ and
$\widetilde{\chi}=0$ near $\pi $,
$$\chi (\theta )T(\theta ) =-2i
\frac{\chi (\theta ) \sin (\theta ) }{\widetilde{\chi} (\theta
)W(\theta )}.$$ By Lemma 2.3 we have $F_0^*W \in \ell ^{1,1+\sigma}
.$ By the argument in Lemma 5.10 \cite{CT} we have
  $F_0^* \left [\frac{W(\theta )}{\sin (\theta )}\right ] \in \ell ^{1,\sigma} .$ Then
$F_0^* \left [ \chi (\theta )T(\theta )\right ] \in \ell ^{1,\sigma}
$ by Wiener's Lemma: case $\sigma =0$ is stated in 11.6 \cite{R};
for $\sigma
>0$ one can provide $\ell ^{1,\sigma}  $ with a structure of
commutative Banach algebra (changing the norm to an equivalent one,
10.2 \cite{R}) and then repeat the argument in 11.6 \cite{R}.

Consider now $A(\theta )= \left ( T(  \theta ) -T(0) \right )  \chi
(\theta )$. We have $F_0^* \left [ A\right ] \in \ell ^{1,\sigma}  $
and $A(0)=A(\pi )=0.$
 We have
 $$\widehat{f}(\nu )=\frac{2i}{\pi } \sum _{\mu    \in \nu +2\Bbb Z+1}
  \frac{1}{\nu -\mu }\widehat{A}(\mu ).$$
We consider
$$  \sum _{\nu \in \Bbb Z}|\widehat{f}(\nu )| \le I+II+III$$
with
$$ \aligned &I=\sum _{\nu \in \Bbb Z} \left |
\sum _{|\mu |\le |\nu |/2, \mu     \in \nu +2\Bbb Z+1}  \frac{
\widehat{A}(\mu ) }{ \nu -\mu   } \right |\, ,\\& II=\sum _{\nu \in
\Bbb Z} \sum _{|\nu |/2\le |\mu |\le 2|\nu | }
\frac{|\widehat{A}(\mu )|}{\langle \nu -\mu\rangle }  \, , \,
III=\sum _{\nu \in \Bbb Z} \sum _{  |\mu |\ge 2|\nu | }
\frac{|\widehat{A}(\mu )|}{\langle \nu -\mu \rangle } .
\endaligned $$
We see immediately that
$$III\lesssim  \|
\widehat{A} \| _{\ell  ^{1,\sigma }} \sum _{\nu \in \Bbb Z}\langle
\nu \rangle ^{-1-\sigma} <\infty .$$ We have
$$II\lesssim \sum _{\mu \in
\Bbb Z} \langle \mu  \rangle ^{ \sigma} | \widehat{A}(\mu ) | \sum
_{| \nu |\le 2|\mu |} \langle \nu  -\mu \rangle ^{-1 }\langle \mu
\rangle ^{ -\sigma}\lesssim \sum _{\mu \in \Bbb Z} \langle \mu
\rangle ^{ \sigma} | \widehat{A}(\mu ) | <\infty .$$ We write
$$\aligned & \sum _{|\mu |\le |\nu |/2, \mu    \in \nu +2\Bbb Z+1}
 \frac{ \widehat{A}(\mu ) }{
\nu -\mu   } =\sum _{|\mu |\le |\nu |/2, \mu    \in \nu +2\Bbb Z+1}
\frac{ \widehat{A}(\mu ) }{ \nu   } +\\& \sum _{|\mu |\le |\nu |/2,
\mu    \in \nu +2\Bbb Z+1}     \frac{\mu }{(\nu -\mu ) \nu } {
\widehat{A}(\mu ) } \ .
\endaligned
$$
Notice $$\aligned &\sum _{\nu \in \Bbb Z}\sum _{|\mu |\le |\nu |/2
}\frac{|\mu \widehat{A}(\mu )| }{\langle \nu -\mu \rangle \langle
\nu \rangle }\lesssim \sum _{\mu \in \Bbb Z} |\mu \widehat{A}(\mu )|
\sum _{|\nu |\ge 2 |\mu | } \langle  \nu \rangle ^{-2}\lesssim \|
\widehat{A}\| _{\ell ^1}<\infty .
\endaligned
$$
The fact that $A(0)=0$ implies $\sum \widehat{A}(\mu ) =0$. The fact
that $A(\pi )=0$ implies $\sum (-1)^\mu \widehat{A}(\mu ) =0$. Hence
$$\sum _{\mu \in 2\Bbb Z}\widehat{A}(\mu )
 =\sum _{\mu \in 2\Bbb Z+1}\widehat{A}(\mu )=0.$$
 This implies that
 $$\sum _{|\mu |\le |\nu |/2,\mu    \in \nu +2\Bbb Z+1}  {
\widehat{A}(\mu ) }= -\sum _{|\mu |> |\nu |/2, \mu    \in \nu +
2\Bbb Z+1} { \widehat{A}(\mu ) }.$$ Then
$$\aligned \sum _{\nu \in \Bbb Z\backslash \{ 0 \}}
  \left | \sum _{|\mu |\le |\nu |/2, \mu    \in \nu + 2\Bbb Z+1}  \frac{
\widehat{A}(\mu ) }{ \nu   }\right |  = \sum _{\nu\in \Bbb
Z\backslash \{ 0 \}}  \left | \sum _{|\mu |> |\nu |/2, \mu    \in
\nu + 2\Bbb Z+1} \frac{ \widehat{A}(\mu ) }{ \nu   }\right |.
\endaligned
$$
This can be bounded with the same argument of $III$. Hence we have
shown $\widehat{f}\in \ell ^1.$

\bigskip
\proclaim{Lemma 3.5} Let $q\in \ell ^{1,1+\sigma }$ with $\sigma
>0$. Then
$U_{112}\in B (L^p,L^p)$  for all $p\in [1,\infty ]$.
\endproclaim
{\it Proof.} The proof is similar to  the previous one. Let $g(\mu ,
\theta )=A(\mu ,\theta ) \text{sign}(\theta )$ with $A(\mu ,\theta
)=\left (n_+(\mu , \theta ) -n_+(\mu ,  0 )\right )\chi (\theta )$.
Set $\widehat{g}(\mu ,\cdot )=F^*[g(\mu , \cdot )]$ and
$\widehat{A}(\mu ,\cdot )=F^*[A(\mu , \cdot )]$. It is enough to
show that there exists $b(\nu )$ in $\ell ^1$ such that
$|\widehat{g}(\mu ,\nu )|\le b(\nu )$ for all $\mu \ge 0$ and all
$\nu\in \Bbb Z$. Notice that $F^*[n_+(\mu , \cdot ) -n_+(\mu ,  0
)](\nu )= \chi _{(0,\infty )}(\nu ) B_{+}(\mu , \nu )$ for $\nu \neq
0$ and $=-n_+(\mu , 0 )$ for $\nu =0$. By  Lemma 2.1 we have
$|B_{+}(\mu , \nu )|\le e^{\gamma (0)}\chi _{(0,\infty )}(\nu )\eta
(\nu )$. Hence $|\widehat{A}(\mu ,\nu )|\le h (\nu )$ for all $\mu
\ge 0$ and all $\nu\in \Bbb Z$, with $h\in \ell ^{1,\sigma }$.

 We have
 $$\widehat{g}(\mu ,\nu )=\frac{2i}{\pi } \sum _{\nu '  -\nu \in 2\Bbb Z+1}
  \frac{1}{\nu -\nu ' }\widehat{A}(\mu ,\nu ')=\frac{2i}{\pi } (I+II+III)$$
with
$$ \aligned &I=
\sum _{|\nu ' |\le |\nu |/2, \nu '    \in \nu +2\Bbb Z+1}  \frac{
\widehat{A}(\mu ,\nu ') }{ \nu -\nu '   }  \, ,\\& II=\ \sum _{|\nu
|/2< |\nu ' |\le 2|\nu |, \nu '    \in \nu +2\Bbb Z+1 } \frac{
\widehat{A}(\mu ,\nu ' ) }{\nu -\nu '  }  \, , \\& III= \sum _{ |\nu
' |> 2|\nu | , \nu '    \in \nu +2\Bbb Z+1} \frac{ \widehat{A}(\mu
,\nu ') }{\nu -\nu ' } .
\endaligned $$
We have
$$|III (\mu , \nu )|\lesssim   \|
h \| _{\ell  ^{1,\sigma }}  \langle \nu \rangle ^{-1-\sigma} .$$ We
have
$$|II  (\mu , \nu )|\lesssim \alpha (\nu ):= \sum _{|\nu
|/2< |\nu ' |\le 2|\nu |  }  \frac{ | h(\nu ' ) |}{    \langle \nu
-\nu ' \rangle   } .$$ We write
$$\aligned & \sum _{|\nu ' |\le |\nu |/2, \nu '    \in \nu +2\Bbb Z+1}
 \frac{ \widehat{A}(\mu ,\nu ') }{
\nu -\nu '   } = I_1+I_2\\& I_1=\frac 1{ \nu   }\sum _{|\nu ' |\le
|\nu |/2, \nu '    \in \nu +2\Bbb Z+1}   \widehat{A}(\mu ,\nu ' ) \,
, \quad I_2= \sum _{|\nu ' |\le |\nu |/2, \nu '    \in \nu +2\Bbb
Z+1} \frac{\nu ' }{(\nu -\nu ' ) \nu } { \widehat{A}(\mu ,\nu ' ) }
.
\endaligned
$$
We have $$ I_1(\mu ,\nu ) =-\frac 1{ \nu   }\sum _{|\nu ' |> |\nu
|/2, \nu '    \in \nu +2\Bbb Z+1}   \widehat{A}(\mu ,\nu ' )$$ and
so
$$|I_1(\mu ,\nu )|\lesssim
\| h\| _{\ell ^{1,\sigma}} \langle \nu \rangle  ^{-1-\sigma }.$$
Finally
$$|I_2(\mu ,\nu )|\lesssim \beta (\nu ):=
\sum _{|\nu ' |\le |\nu |/2,  } \frac{\langle \nu '\rangle  }{
\langle \nu -\nu ' \rangle \langle  \nu \rangle } { h(\nu ' ) }
$$
Then there is a function $b(\nu )$ in $\ell ^1$ such that
$|\widehat{g}(\mu ,\nu )|\le b(\nu )$ of the form $b(\nu )= C
(\alpha (\nu )+\beta (\nu )+ \langle  \nu \rangle ^{-1-\sigma})$.

\bigskip
By repeating the previous arguments one has:

\proclaim{Lemma 3.6} For $q\in \ell ^{1,2+\sigma}$ with $\sigma
>0$ the operator $W$ extends into a bounded operator in $\ell ^p$
for $p=1,\infty$ when operators (3.1)--(3.4) are bounded. Here (3.1)
has been defined above while (3.2)--(3.4) are defined as follows,
for $\chi +\chi _1$ a smooth partition of unity in $\Bbb T$ with
$\chi =1$ near 0 and  $\chi =0$ near $\pi$:
$$  \aligned & V_{2}f(\mu ) = \chi _{[0,\infty )}(\mu )\left (
T(\pi )-1\right ) m_+(\mu , 0 ) \int _{-\pi }^\pi
  e^{  i\mu \theta } \text{sign} (\theta )  \chi _1(\theta )
    F_0[f] ( \theta )
     \\&  - \chi _{[0,\infty )}(\mu ) R_+(  \pi  )m_+(\mu ,
0 ) \int _{-\pi }^\pi
  e^{ - i\mu \theta } \text{sign} (\theta )\chi _1(\theta ) F_0[f] (  \theta )
d\theta   .\endaligned \tag 3.2$$

$$  \aligned & V_{3}f(\mu ) = \chi _{(-\infty ,0)}(\mu )\left (
1-T(0)\right ) m_-(\mu , 0 ) \int _{-\pi }^\pi
  e^{  i\mu \theta } \text{sign} (\theta )  \chi  (\theta )
    F_0[f] ( \theta )
     \\&  + \chi _{(-\infty ,0)}(\mu ) R_-(  0 )m_-(\mu ,
0 ) \int _{-\pi }^\pi
  e^{  -i\mu \theta } \text{sign} (\theta )\chi  (\theta ) F_0[f] (  \theta )
d\theta   .\endaligned \tag 3.3$$

$$  \aligned & V_{4}f(\mu ) = \chi _{(-\infty ,0)}(\mu )\left (
1-T(0) \right ) m_-(\mu , 0 ) \int _{-\pi }^\pi
  e^{  i\mu \theta } \text{sign} (\theta )  \chi _1(\theta )
    F_0[f] ( \theta )
     \\&  + \chi _{(-\infty ,0)}(\mu ) R_-(  0 )m_-(\mu ,
0 ) \int _{-\pi }^\pi
  e^{  -i\mu \theta } \text{sign} (\theta )\chi _1(\theta ) F_0[f] (  \theta )
d\theta   .\endaligned \tag 3.4$$
\endproclaim

We have: \proclaim{Lemma 3.7} $W\in B(\ell ^p,\ell ^p)$ for
$p=1,\infty$ exactly when $ T(0)= T(\pi )=1.$
\endproclaim
{\it Proof.} If $ T(0)= T(\pi )=1 $ we have $V_j=0$ for all $j$.
Then $W\in B(\ell ^p,\ell ^p)$ for $p=1,\infty$. Viceversa
 $W\in B(\ell ^1,\ell ^1)$ implies  $V_j\in B(\ell ^1,\ell ^1)$ for
all $j$. If $V_3\in B(\ell ^1,\ell ^1)$ then, proceeding as in Lemma
3.3,
$$ 1-T(0)-R_-(0)= 1-T(0)+R_+(0)= 0.$$
 This together with
(1) in Lemma 3.3 implies $ T(0)=  1.$ The implication $ T(\pi )= 1 $
is obtained similarly.

\bigskip

\head \S A Appendix: finite number of eigenvalues \endhead

We will prove:
 \proclaim{Lemma A.1} If $q\in
\ell ^{1,1}$ the total number of eigenvalues of $H$      is $\le 4+
\|\nu q (\nu )\| _{\ell ^1}.$
\endproclaim
Let $q _-(\nu )=\min ( 0, q(\nu ))$. We recall that if we have
$(-\Delta +q)u=\lambda u$, then if we define $v$ by $v(\nu
)=(-1)^{\nu} u(\nu )$ we have $(-\Delta -q)v=(4-\lambda )v$. Hence
Lemma 6.1 is a consequence of:

 \proclaim{Lemma A.2} If $q\in
\ell ^{1,1}$ the total number of eigenvalues of $H$  inside
$(-\infty , 0)$   is $\le 2+ \|\nu q_- (\nu )\| _{\ell ^1}.$
\endproclaim
 {\it Proof.} For $\lambda \le 0$ we set
  $u(\nu ,\lambda )=f _+(\nu , \theta )$, where $\lambda
  =2(1-\cos (\theta )).$ Notice that   $u(\nu ,\lambda )\in \Bbb
  R$. We denote by $X(\lambda )$   the set of those
  $\nu $ such that either $u(\nu ,\lambda )=0$ or
$u(\nu ,\lambda )u(\nu +1 ,\lambda ) <0$. We denote by $N(\lambda )$
the cardinality of $X(\lambda )$.  Notice that by the min-max
principle the operator $\widetilde{H}=-\Delta - q _-$   has at least
as many negative eigenvalues as $H$. So, to prove our Lemma 6.2 it
is not restrictive to assume $q(\nu )=q_-(\nu )=-| q (\nu )| $ for
all $\nu$ in Lemma  A.3 below. We have: \proclaim{Lemma 6.3} We have
$N(0)\le 2+ \|\nu q _-(\nu )\| _{\ell ^1}.$
\endproclaim
{\it Proof.} We assume $N(0)>1$. Let $\nu _0,\nu _1\in X (0)$ be two
consecutive elements, with $\nu _0<\nu _1$. For
 $u(\nu )=u(\nu , 0)$ we
have $$u(\nu  )= u (\nu _0)  +(u (\nu _0+1)-u (\nu _0)) (\nu -\nu
_0)-\sum _{j=\nu _0}^{\nu -1}(j -\nu _0)|q(j)| u(j  ). $$ It is not
restrictive below to assume $A:= u (\nu _0+1)-u (\nu _0) >0$. Then
$u (\nu _1+1)<0$ or $u (\nu _1 )=0$. In the first case,   we have
$$0>u (\nu _0+1)-u (\nu _1+1)=A (\nu _1 -\nu _0)\left (1 -
\sum _{j=\nu _0}^{\nu _1}(j -\nu _0)|q(j)|  \right ) .
$$
This implies
$$\sum _{j=\nu _0+1}^{\nu _1}(j -\nu _0)|q(j)|\ge 1.
\text{ By a similar argument } \sum _{j=\nu _0 }^{\nu _1-1}( \nu
_1-j)|q(j)|\ge 1 . \tag 1$$ (1) holds also if $u (\nu _1 )=0$. So
for $ \nu _0<\nu_1<...<\nu _n$ consecutive elements in $X(0)$,
$$\text{we have }\sum _{j=\nu _0+1}^{\nu _n}(j -\nu _0)|q(j)|\ge n
 \text{ and   }
\sum _{j=\nu _0 }^{\nu _n-1}( \nu _n-j)|q(j)|\ge n.$$ Then $q\in
\ell ^{1,1}$ implies $N(0)<\infty$. If  $X(0)$ is formed by
$$ \nu _0<...<\nu _{n}(<0\le ) \mu _0<...<\mu _{m}$$
then
$$n\le \sum _{j=\nu _0 }^{\nu _n-1}( \nu _n-j)|q(j)| \le
 \sum _{j=\nu _0 }^{\nu _n-1}|j||q(j)| $$ and
 $$m\le \sum _{j=\mu _0+1}^{\mu _m}(j -\mu _0)|q(j)|\le
  \sum _{j=\mu _0+1}^{\mu _m} |j||q(j)| .$$
 So $n+m\le \| \nu q(\nu ) \| _{\ell ^1}$. Then
  $N(0)\le 2+\| \nu q(\nu ) \| _{\ell ^1}$. This yields Lemma 6.2.

\bigskip
  Notice that
$$\langle H u,u\rangle =\sum _{\nu \in \Bbb Z}|u(\nu +1)-u(\nu )|^2+
\sum _{\nu \in \Bbb Z}q(\nu )|u(\nu )|^2.$$ If $H$ has negative
eigenvalues, there is a minimal one $\lambda _0$. Then we have
$u(\nu ,\lambda _0)=|u(\nu ,\lambda _0)|>0$ for all $\nu$    by the
min-max  principle and by the fact that $u(\nu ,\lambda _0)=e^{i\nu
\theta}m_+(\nu ,\theta _0)$ where $m_+(\nu ,\theta )\to 1$ for $
|\nu |\nearrow   \infty$ by (1) Lemma 5.1 \cite{CT}. Notice that by
this argument it is easy to conclude that $N(\lambda )<\infty $ for
any $\lambda < 0$.

 Next we
have the following discrete version of the Sturm oscillation
theorem, see Lemma 4.4 \cite{T}.
   \proclaim{Lemma A.4} $N(\lambda )$ is increasing  for $\lambda
   \le 0$.
\endproclaim
Lemmas A.4 and A.3 yield Lemma A.2.


\Refs\widestnumber\key{1997shire}

\ref\key{C} \by S.Cuccagna \paper On  asymptotic  stability in 3D of
kinks for the $\phi ^4$ model \jour  Trans. Amer. Math. Soc.
 \vol 360
\yr 2008 \pages   2581-2614
\endref


\ref\key{CT} \by S.Cuccagna, M.Tarulli  \paper On asymptotic
stability   of  standing waves of discrete
 Schr\"odin- ger equation  in $\Bbb Z$
\jour      \vol   \yr   \pages
\endref

\ref\key{DF} \by P.D'Ancona, L.Fanelli  \paper  $L^{ p}$ boundedness
 of the wave operator  for the one dimensional   Schr\"odinger
 operator
\jour   Comm. Math. Phys.  \vol  268  \yr 2006  \pages 415--438
\endref





\ref\key{DT} \by P.Deift, E.Trubowitz \paper Inverse scattering on
the line\jour Comm. Pure Appl. Math. \vol 32 \yr 1979 \pages
121--251
\endref


\ref\key{GY} \by Galtabiar, K.Yajima \paper $L^p$ boundedness of
wave operators for one dimensional Schr\"o-dinger operators \jour J.
Math. Sci. Univ. Tokio \vol 7 \yr 2000 \pages 221 -- 240
\endref



\ref\key{KKK} \by A.Komech, E.Kopylova, M.Kunze \paper Dispersive
estimates for 1D discrete Schr\"odinger and Klein Gordon
equations\jour Appl. Mat. \vol 85 \yr 2006 \pages 1487--1508
\endref



\ref\key{GSc}\by  M.Goldberg, W.Schlag \paper Dispersive estimates
 for Schr\"odinger   operators in dimensions one and three
\jour  Comm. Math. Phys. \vol 251 \yr 2004 \pages 157--178
\endref





\ref\key{PS} \by D.E.Pelinovsky,A. Stefanov  \paper On the spectral
theory and dispersive estimates for a discrete Schr\"odinger
equation in one dimension \paperinfo
\endref






\ref\key{RS} \by M.Reed, B.Simon \book Methods of mathematical
Physics\publ Academic Press \publaddr   San Diego \yr 1979
\endref


\ref\key{R}\by W.Rudin \book Functional Analysis \bookinfo Higher
Math. Series \publ McGraw-Hill \yr 1973
\endref


\ref\key{T}\by G.Teschl \book Jacobi Operators and Completely
Integrable Nonlinear Lattices \bookinfo Mathematical Surveys  and
Monographs \publ AMS \yr 2000
\endref






\ref\key{SK} \by A. Stefanov, P.G.Kevrekidis  \paper Asymptotic
behaviour of small solutions for the discrete nonlinear
Schr\"odinger and Klein--Gordon equations
 \jour Nonlinearity \vol 18 \yr 2005 \pages
1841--1857
\endref

\ref\key{W1} \by R. Weder \paper The $W^{k,p}$ continuity of the
Schr\"odinger
  wave operators on the line
\jour   Comm. Math. Phys. \vol 208 \yr 1999 \pages 507--520
\endref

\ref\key{W2} \bysame \paper $L^p\to L^{p^\prime}$ estimates for
 the Schr\"odinger equation
   on the line and inverse
scattering for the nonlinear Schr\"odinger equation with a potential
\jour   J. Funct. Anal. \vol 170 \yr 2000 \pages 37--68
\endref






\endRefs
\enddocument